\documentclass[12pt,a4paper]{article}
\usepackage{amssymb, amsmath, amsthm, amscd, graphicx, amsfonts}
\newtheorem{theo}{Theorem}
\begin{document}
\begin{center}
{\textsf{\large On Long Range Percolation With Heavy Tails}}\\
\vspace{3mm}
S. Friedli\footnote{CBPF, Rua Dr. Xavier Sigaud, 150, Urca, CEP 22290-180, 
Rio de Janeiro, and 
IMPA, Estrada Dona Castorina, 110 Jardim Bot\^anico CEP 22460-320, 
Rio de Janeiro.}, B.N.B. de Lima\footnote{UFMG, Av. Ant\^onio Carlos 6627 
C.P. 702  
Belo Horizonte 30123-970.}, V. Sidoravicius\footnote{IMPA, Estrada Dona 
Castorina, 110 Jardim Bot\^anico 
CEP 22460-320, Rio de Janeiro.}
\end{center}
\begin{abstract}
Consider independent long range percolation on $\mathbf{Z}^2$,
where horizontal and vertical edges of length $n$ are open with
probability $p_n$. We show that if
 $\limsup_{n\to\infty}p_n>0,$ then there
exists an integer $N$ such that $P_N(0\leftrightarrow \infty)>0$,
where $P_N$ is the truncated measure obtained by taking
$p_{N,n}=p_n$ for  $n \leq N$ and $p_{N,n}=0$ for all $n> N$.
\end{abstract}



\vspace{0.3cm}

On the graph ${\cal G}=({\cal V},{\cal E})$ with  ${\cal
V}=\mathbf{Z}^2$, and  ${\cal E}= \{ \langle x,y \rangle \subset
{\cal V} \times {\cal V}: \; |x_1 - y_1| = 0, \; {\text {or}} \;
|x_2 - y_2|=0\}$, consider a  long range percolation process
$(\Omega,{\cal F},P)$, where $\Omega=\{0,1\}^{\cal E}$,
$P=\prod_{\langle x,y \rangle \in{\cal E}}\mu_{\langle x,y
\rangle}$, and $\mu_{\langle x,y \rangle}\{\omega_{\langle x,y
\rangle} = 1\}=p_{|x-y|} \in [0,1]$ is a Bernoulli measure,
independent of the state of other edges. Given a sequence
$(p_n)_{n\in\mathbf{N}}$ and $N\in\mathbf{N}$, we
define a truncated sequence $(p_{N,n})_{n\in\mathbf{N}}$ by
\begin{equation}\label{trunc}
p_{N,n}=
\begin{cases}
p_n&\text{ if }n\leq N\,,\\
0 &\text{ if } n>N,
\end{cases}
\end{equation}
and a truncated percolation process by taking
$P_N=\prod_{\langle x,y \rangle\in\cal{E}}\mu_{N,\langle x,y
\rangle}$, where $\mu_{N,\langle x,y \rangle}\{\omega_{\langle x,y
\rangle} = 1\}=p_{N,|x-y|}$.

In this note we adress the following question: given a
sequence $(p_n)_{n\in\mathbf{N}}$ for which 
$P(0\leftrightarrow\infty)>0$, does
there always exists some large enough $N$ such that
$P_N(0\leftrightarrow\infty)>0$?
In other words, given a
system with infinite range translation invariant interactions which exhibits a
phase transition, we ask if the infiniteness of the range is really
crucial for this transition to occur. 
It is known, for instance, that infinite range is essential
in one dimensional systems
(cf. \cite{FS}, \cite{NS}), but
it is believed that in dimensions $d\geq 2$,
occurrence  (if so) of
phase transitions for translation invariant interactions is  {\it always} determined by a bounded part
of the interaction (excluding cases when interactions are of intrinsically one dimensional structure).
Returning to the percolation case, rapid (say, exponential) decay or summability of the $p_n$'s indicates that long range
connections may not be necessary for the existence of an infinite cluster. 
This is the setup of \cite{MS} and partially \cite{Be}. On the other hand, 
heavy tail interactions are still poorly
understood, and the only existing studies, \cite{SSV} and \cite{Be}, 
rely heavily on asymptotic monotonicity assumptions which are a key
ingredient for the use of rather laborious coarse-graining techniques.
Although the approach we present here relies on deep and highly nontrivial
facts (\cite{GM}, \cite{K}),
it leads to a much shorter (not to say
elementary) proof, which is
less sensitive to the geometry of the interactions, and
allows to consider a rather general class of systems with connections of
irregular, in particular lacunary structure.

\begin{theo}\label{PTT1}  $(d\geq 2)$ If
$\limsup_{n\to\infty}p_n>0,$ then
\begin{equation}\label{PT1.1}
P_N(0\leftrightarrow\infty)>0
\end{equation}
for some large enough $N$.
\end{theo}

\noindent {\bf Proof.} It suffices to consider $d=2$. Define $\epsilon>0$ by
$2\epsilon =\limsup_{n\to\infty}p_n\, >0$. By \cite{K},
Theorem 1 p. 220, there exists
$d_{\epsilon}$ such that
\begin{equation}\label{PT5}
p_c(\mathbf{Z}^{d_\epsilon}) < \epsilon/2\,,
\end{equation}
and by \cite{GM}, Theorem A p. 447, we can find $K_\epsilon$ such that
\begin{equation}\label{PT6}
p_c\big(\{1,2,\dots,K_\epsilon\}^{d_\epsilon-2}\times
\mathbf{Z}^2\big) < p_c(\mathbf{Z}^{d_\epsilon}) + \epsilon/2 <
\epsilon\,.
\end{equation}
\noindent Let $n_0=0$. 
For $j \in \{1,2, \dots , d_\epsilon -1\}$, define
recursively
$$n_j
= \min \{\ell
> (K_\epsilon + 1)n_{j-1} : \; p_\ell \geq \epsilon\}\,.
$$
For $x \in \mathbf{Z}^2$ and $B \subseteq \mathbf{Z}^2$, define
$\mathbf{T}_{x}B = \{ z+x, \; z\in B\}$. Set $B_0 = \{(0,0)\}$.
For $j \in \{1,2, \dots , d_\epsilon -2\}$, define $B_j
=\cup_{m=0}^{K_\epsilon -1} \mathbf{T}_{m(n_j,0)}B_{j-1}$. Then, let
$${\cal V}_{d_\epsilon -1} = \bigcup_{(k,m) \in
\mathbf{Z}^2}\mathbf{T}_{(k n_{d_\epsilon-1},m
n_1)}B_{d_\epsilon-2}\,,
$$
and
\begin{align}
{\cal E}_{d_\epsilon -1} = \{&\langle x, y\rangle, \;
x,y \in {\cal V}_{d_\epsilon -1} :
 |x_1 -y_1| = n_j\, {\text {for some }} 
1 \le j \le d_\epsilon
-1\nonumber\\
& {\text {and}} \; x_2= y_2,\, 
 {\text {or}}\,\,  x_1=y_1 \; {\text {and}}\;
|x_2 - y_2| = n_1 \}\,.\nonumber
\end{align}

\noindent It 
is straightforward that the graph 
${\cal G}_{d_\epsilon -1}= \big( {\cal
V}_{d_\epsilon -1},\, {\cal E}_{d_\epsilon -1} \big)$
is isomorphic to
$\{1,2,\dots,K_\epsilon\}^{d_\epsilon-2}\times \mathbf{Z}^2$.
Moreover, by our choice of $n_j, \; 1 \le j \le d_\epsilon -1$, we
have that each edge of ${\cal G}_{d_\epsilon -1}$ is open with
probability at least $\epsilon$, and using \eqref{PT6} we get
\eqref{PT1.1} with $N=n_{d_\epsilon -1}$. \qed
\medskip

\noindent {\bf Remark.} For further applications of this method to 
percolation and interacting spin systems  see \cite{FN}.

\smallskip

\noindent {\bf Acknowledgments.} S.F. is supported by the Fonds National 
Suisse pour la Recherche 
Scientifique, B.N.B.L. and V.S. are partially supported by 
CNPq and FAPERJ. F.S. wishes to thank CBPF and IMPA for hospitality 
and support, B.N.B.L. wishes to 
thank IMPA for hospitality and financial support during multiple visits.

\end{document}